\documentclass[12pt]{article}
\usepackage{graphicx,a4}
\usepackage{amsfonts,amsmath,amssymb,mathrsfs,amsthm,epsfig,color}
\definecolor{Red}{cmyk}{0,1,1,0}

\definecolor{verde}{cmyk}{1,0,1,0}

\definecolor{azul}{cmyk}{1,1,0,0}

\newtheorem{theorem}{Theorem}

\newtheorem{lemma}[theorem]{Lemma}
\newtheorem{definition}{Definition}

\newtheorem{remark}[definition]{Remark}

\def\ds{\displaystyle}

\let\a=\alpha \let\b=\beta   
 \let\g=\gamma      
       
 \let\s=\sigma  
  
   \let\L=\Lambda

\global\newcount\numsec\global\newcount\numfor
\gdef\profonditastruttura{\dp\strutbox}
\def\senondefinito#1{\expandafter\ifx\csname#1\endcsname\relax}
\def\SIA #1,#2,#3 {\senondefinito{#1#2}
\expandafter\xdef\csname #1#2\endcsname{#3} \else
\write16{???? il simbolo #2 e' gia' stato definito !!!!} \fi}
\def\etichetta(#1){(\veroparagrafo.\veraformula)
\SIA e,#1,(\veroparagrafo.\veraformula)
 \global\advance\numfor by 1
 \write16{ EQ \equ(#1) ha simbolo #1 }}
\def\etichettaa(#1){(A\veroparagrafo.\veraformula)
 \SIA e,#1,(A\veroparagrafo.\veraformula)
 \global\advance\numfor by 1\write16{ EQ \equ(#1) ha simbolo #1 }}
\def\BOZZA{\def\alato(##1){
 {\vtop to \profonditastruttura{\baselineskip
 \profonditastruttura\vss
 \rlap{\kern-\hsize\kern-1.2truecm{$\scriptstyle##1$}}}}}}
\def\alato(#1){}
\def\veroparagrafo{\number\numsec}\def\veraformula{\number\numfor}
\def\Eq(#1){\eqno{\etichetta(#1)\alato(#1)}}
\def\eq(#1){\etichetta(#1)\alato(#1)}
\def\Eqa(#1){\eqno{\etichettaa(#1)\alato(#1)}}
\def\eqa(#1){\etichettaa(#1)\alato(#1)}
\def\equ(#1){\senondefinito{e#1}$\clubsuit$#1\else\csname e#1\endcsname\fi}

\def\0{\emptyset}

\def\\{\noindent}

\def\tende#1{\vtop{\ialign{##\crcr\rightarrowfill\crcr
              \noalign{\kern-1pt\nointerlineskip}
              \hskip3.pt${\scriptstyle #1}$\hskip3.pt\crcr}}}
\def\otto{{\kern-1.truept\leftarrow\kern-5.truept\to\kern-1.truept}}

\def\1{\rlap{\mbox{\small\rm 1}}\kern.15em 1}

\def\buildd#1#2{\mathrel{\mathop{\kern 0pt#1}\limits_{#2}}}

\def\Z{\mathbb{Z}}

\BOZZA

\def\be{\begin{equation}}
\def\ee{\end{equation}}

\newcommand{\eop}{\nopagebreak\hfill\fbox{ }}

\begin{document}
\begin{center}
\vspace*{0.7 cm}
{\large\bf  Diffusive-Ballistic Transition
in Random Polymers \vspace*{.2cm}\\with Drift and Repulsive Long-Range Interactions}
\vspace*{0.5cm}\\
    L.M. Cioletti$^{\mbox{\footnotesize{a},}}$
    \footnote{Corresponding author : leandro.mat@gmail.com
    				\\ 
			  	\hspace*{0.55cm} Co-authors : changdorea@unb.br ; simone@ufg.br
			  	\\
				\hspace*{0.55cm} Research partially supported by CNPq and FEMAT.
	},
    C.C.Y. Dorea$^{\mbox{\footnotesize{a}}}$ and S. Vasconcelos da Silva$^{\mbox{\footnotesize{b}}}$
\\
 \vskip 3mm
$\phantom{.}^{ \mbox{\footnotesize{a}}}$ 
Universidade de Brasilia, 70910-900 Brasilia-DF, Brazil
\\
$\phantom{.}^{ \mbox{\footnotesize{b}}}$ 
Universidade Federal de Goias, 74001-970 Goiania-GO, Brazil
\end{center}
\vskip 5mm
 \begin{abstract}
{\footnotesize
In this note phase transition issues are addressed 
for random polymers on $\Z^2$ with long-range self-repulsive interactions. 
It is shown that, in the absence of drift and with power law interactions, 
the polymer exhibits transition from diffusive to a ballistic behavior. 
When non-null drifts are added and positive translation invariant 
interactions are considered, the polymer presents a ballistic behavior. 
Our results complement some previous studies on the matter and  
we also derive a Central Limit Theorem for the model.
}
\end{abstract}
\vskip 1mm
{\footnotesize 
MSC: 82B20, 82B41, 82B26.
\newline
Keywords: self-repelling random polymers; 
Ising model; 
long-range interactions;
diffusive-ballistic phase transition; CLT.
}

\section{Introduction}
\indent 

Random polymers can be modelled as connected subsets of $\Z^2$. More precisely, a $N$-th step polymer $S$
is an element of $\mathbb{W}_N$ given by  
\[
\mathbb{W}_N:= \{S=(S_0,S_1,\dots, S_N) : S_i\in \Z^2, S_0=0 \ 
\text{and}\ \|S_{i+1}-S_i\|=1\},
\]
being $\|\cdot\|$ the $\ell^1$ norm. Under Gibbs measure setting at inverse temperature $\beta>0$ and Hamiltonian $\mathbb{H}_N$ we can write the probability
\be\label{model}
\mathbb{P}_{N}^{\b,h}(S)
=
\frac{\exp 
\left[-\b\mathbb{H}_{N}(S)\right]}{\mathbb{Z}^{\b}_N(h)}\ ,
\qquad
\mathbb{H}_{N}(S) = - \!\!\!\!
\sum_{{ 1\leq i<j\leq N}}\!\!\!
V_{ij}\langle X_i,X_j\rangle
+ 
\langle h,S_N \rangle,
\ee
where $X_i=S_{i}-S_{i-1}$ stands for the $i$-th random step, $V_{ij}$ are the prescribed interactions, $\langle\cdot,\cdot\rangle$ denotes the usual inner product, $h\in\mathbb{R}^ 2$ is the fixed drift vector and $\mathbb{Z}_N^{\b}(h)$ is the partition function.

Caracciolo et al. (\cite{cpp}, 1993) 
introduced a self-repelling random polymer model 
with Hamiltonian on $\mathbb{W}_N$ given by
$
\overline{H}_{N}(S)=g_0\sum_{0\leq i<j\leq N}V_{ij}\delta_{S_i,S_j},
$
where $g_0> 0$ and the interactions $V_{ij}=|i-j|^{-\alpha}$. 
Their model interpolates between the lattice Edwards model ($\a=0$) and ordinary SRW ($\a=\infty$). 
Moreover, it was conjectured that for dimension $1\leq d \leq 4$ there exists a strictly positive exponent 
$\gamma=\gamma(d,\alpha)$ such that the mean square end-to-end distance satisfies the asymptotics
$$
\ds\mathbb{E}_{\overline{P}_N}[\|S_N\|^2]=\sum_{S\in\mathbb{W}_N}\|S_N\|^2 \overline{P}_N(S)\sim cN^{\gamma},
$$
where the Gibbs measure $\overline{P}_N$ 
is given by the Hamiltonian $\overline{H}_N$.
In \cite{arb}, the Hamiltonian 
$
\tilde{H}_N(S)=-\sum_{0\leq i<j\leq N}{|i-j|^{-\a}}\|S_i-S_j\|^2$, 
	where $ 3<\alpha\leq 4$, 
was considered. They proved the existence of positive constants 
$\b_1$ and $\b_2$ that led to phase transition from diffusive regime ($\b<\b_1$) to  a ballistic one ($\b>\b_2$). However, it was left unknown what undergoes when $\b\in [\b_1, \b_2]$. As usual, the different diffusive regimes are classified according to the asymptotic behavior of the mean square displacement and for our model  \eqref{model} it reduces in determining $\g>0$ for which the following limit exists, is positive and finite
\be
\label{quadrado}
\lim_{N\to \infty} \frac{1}{N^\g} 
{\mathbb{E}_{\mathbb{P}_{N}^{\b,h}}[\|S_N\|^2]}= \lim_{N\to \infty} \frac{1}{N^\g}{\sum_{S \in \mathbb{W}_N}\|S_N\|^2\ \mathbb{P}_{N}^{\b,h}(S)}.
\ee
We say that the polymer model 
 is  {\it diffusive} if $\g=1$, {\it
superdiffusive} if $1<\g<2$ and {\it ballistic} if $\g=2$.

Our main motivation is to build a self-repelling random polymer model 
for which we can derive a genuine diffusive-ballistic phase transition, i.e. 
the existence of a unique positive constant $\beta_c$ separating the model into two regimes.
In this note, assuming zero drift and $V_{ij}=|i-j|^{-\a}$ with $1<\alpha\leq 2$, 
we prove (Theorem \ref{teo:phasetransition})
that there exists a unique positive number $\b_c$ (the critical temperature of a related 
one dimensional Ising model) such that the model is diffusive for $\b<\b_c$ and ballistic for $\b>\b_c$. 
On the other hand, considering non-null drift 
and positive, translation invariant 
and regular interactions, 
we conclude from Theorem \ref{teo:h1h2-positivo}
that for all $\beta\in (0,\infty)$ the model
is ballistic. 

The Lemma \ref{lem-esp} is an essential tool in this work. 
Its proof is similar in spirit to the one introduced for the Potts model  
by M. Suzuki \cite{moto} in 1967. It consists in decoupling the steps of the 
polymer as two independent Ising random variables. 
The background idea is the 
same applied when looking at SRW in lattice $\mathbb{Z}^2$ 
as two independents SRW's on $\mathbb{Z}$.
 
In 1983 Newman \cite{Newman-CLT} proved a CLT for block random variables  
satisfying the FKG inequalities under finite susceptibility hypothesis. In Section 3 we investigate the validity of the CLT for our model. 
Here assuming non zero drifts and consequently infinite susceptibility, 
we prove (Theorem \ref{teo:clt}) a CLT for the projections of suitably normalized displacements. 
This is obtained by using 
both the Lee-Yang circle theorem and the $C^2$-regularity condition from Wu Liming \cite{WL}. 
 
\section{Mean Square Displacement and Phase Transition}
\indent

For the volume $\Lambda_N = \{1,2,\dots ,N\}$ consider the one dimensional Ising model  with free boundary conditions 
defined by the Hamiltonian 
$$
H_{\Lambda_N}(\s)=-\sum_{1\le i<j\le N}V_{ij}\s_i\s_j-\sum_{i=1}^N h\sigma_i,
$$
where $\sigma=(\sigma_1,\ldots,\sigma_N)\in \{-1,1\}^N:= \Sigma_N$, 
$V_{ij}\in \mathbb{R}$ are the coupling constants 
and $h\in\mathbb{R}$ is an external field.
To simplify notation for a given a real-valued function $f:\Sigma_N\to\mathbb{R}$ write
$$
\left\langle f\right\rangle_{\Lambda_N}^{\b,h}
= \mathbb{E}_{P_{{\Lambda_N}}^{\b,h}}[f] \quad  \text{with} \quad
P_{{\Lambda_N}}^{\b,h}(\s)= 
\frac{1}{Z_{\L_N}^{\b}(h)}
\exp\left(-\b H_{\Lambda_N}(\s)\right)
$$
where $Z_{\L_N}^\b(h)$ is the partition function.

\begin{lemma}\label{lem-esp}
For $e_1=(1,0)$ and $e_2=(0,1)$ define $h_1=\langle h,e_1-e_2\rangle$ and $h_2=\langle h,e_1+e_2\rangle$. 
Then 
\[
\mathbb{E}_{\mathbb{P}_{N}^{\b,h}}[\|S_N\|^2]
=
\frac{1}{2}\sum_{i,j=1}^N
\Big[
 \left\langle \s_i\s_j\right\rangle_{\Lambda_N}^{\frac{\b}{2},h_1}
+
\left\langle \s_i\s_j\right\rangle_{\Lambda_N}^{\frac{\b}{2},h_2}
\Big].
\]
\end{lemma}
\noindent{\bf Proof.} 
The proof follows closely the ideas from \cite{bps,arb}.
Let $T:\mathbb{R}^2\to\mathbb{R}^2$ be the rotation
$
TX_i=(\s_i e_1+\tilde\s_i e_2) /\sqrt 2,
$
with $\s_i,\tilde\s_i\in\{-1,1\}$. A simple computation shows that 
\vspace{-0.4cm}
$$
\mathbb{P}_{N}^{\b,h}(S)=
P_{{\Lambda_N}}^{\frac{\b}{2},h_1}(\s)P_{{\Lambda_N}}^{\frac{\b}{2},h_2}(\tilde\s)\quad  \text{and} \quad 
\|S_N\|^2
= 
\frac{1}{2}\sum_{i,j=1}^N
(\s_i\s_j+\tilde\s_i\tilde\s_j).
\vspace{-0.4cm}
$$
\eop

\begin{theorem}\label{teo:h1h2-positivo}
Suppose that $V_{ij}$ is positive and  tranlation invariant, i.e. $V_{ij}=V(|i-j|)>0$ for all $i\neq j$. 
If $h\in \mathbb{R}^ 2$ is such that $h_1$ and $h_2$   
satisfy $h_1h_2>0$ and $\beta>0$, 
then for some constant $C(\beta,h)> 0$ we have  
$$
C(\beta,h) 
\leq 
\frac{\mathbb{E}_{\mathbb{P}_{N}^{\b,h}}[\|S_N\|^2]}{N^2}.
$$
\end{theorem}

\noindent{\bf Proof.} Let $k\in\mathbb{R}$. Since $V_{ij}>0$  for $i\neq j$  
we get from the second Griffiths inequality, 
\[
\left\langle \s_i\right\rangle_{\Lambda_N}^{\b,k,nn}
\left\langle \s_j\right\rangle_{\Lambda_N}^{\b,k,nn}
\leq
\left\langle \s_i\s_j\right\rangle_{\Lambda_N}^{\b,k,nn}
\leq
\left\langle \s_i\s_j\right\rangle_{\Lambda_N}^{\b,k}
\]
where the left hand side expected values are
taken with respect to the Gibbs measure of the nearest neighbours Ising model on $\Lambda_N$
with free boundary conditions and Hamiltonian given by
$
H_{\Lambda_N}(\sigma)=-\sum _{n=1}^{N-1}V(1)\sigma_n\sigma_{n+1}-k\sum_{n=1}^N \sigma_n .
$ 
Using monotonicity with respect to the volume and classical transfer matrix computation (see \cite[p. 107]{ellis})
we have for any $i\in \Lambda_N$ 
$$
\left\langle \s_i\right\rangle_{\Lambda_N}^{\b,k,nn}\underset{N\to\infty}{\longrightarrow}
\frac{\sinh \beta k}{\sqrt{\sinh^2(\beta k)+e^{-4\beta V(1)}}}.
$$ 
\eop 

\begin{theorem}\label{teo:phasetransition}

Let $1< \alpha\leq 2$, $h=0$ and $V_{ij}=|i-j|^{-\alpha}$ for $i\neq j$.  
Then there exist a constant $\beta_c\in(0,\infty)$ and positive numbers $m_*(\b)$ and $K(\beta)$ 
such that
\be\label{ballis}
\frac{1}{2}m_*^2(\b)\leq \frac{1}{N^2}\mathbb{E}_{\mathbb{P}_{N}^{2\b,0}}[\|S_N\|^ 2]\leq 1, 
\qquad\text{if}\quad \beta > \b_c
\ee
and
\be\label{diffu}
1\leq \frac{1}{N}\mathbb{E}_{\mathbb{P}_{N}^{2\b,0}}[\|S_N\|^ 2]\leq K(\beta) , \qquad\text{if}\quad 0<\beta<\beta_c.
\ee
\end{theorem}

\noindent{\bf Proof.} 
For $1<\a\le 2$ the existence of a critical $\b_c\in(0,\infty)$ for the long range Ising model with coupling 
$V_{ij}$ is shown in  \cite{d, fs}. In this case, we have spontaneous magnetization $m_*(\b)>0$  
for all $\b>\b_c$ and the two-point function
with free boundary condition satisfies (cf. \cite{l}), 
$$
\langle \s_i\s_j\rangle^{\beta,0} 
=
\frac{1}{2}
\Big[ \langle \s_i\s_j\rangle^{\beta,0,+} + \langle \s_i\s_j\rangle^{\beta,0,-}\Big]
\geq m_*^2(\beta)\geq m_*^2(\beta_c).
$$
Using the same type of arguments as in Theorem \ref{teo:h1h2-positivo} we have for large $N$ 
\[
\mathbb{E}_{\mathbb{P}_{N}^{2\b,0}}[\|S_N\|^2]
=
\sum_{i,j=1}^{N}\left\langle \s_i\s_j\right\rangle_{\Lambda_N}^{\b,0}
\ge\frac{1}{2}\sum_{i,j=1}^{N}\left\langle \s_i\s_j\right\rangle^{\b,0}\ge \frac{1}{2}m^2_*(\beta) N^2.
\]

To prove (\ref{diffu}) one needs lower and upper bounds for 
$\left\langle \s_i\s_j\right\rangle_{{\Lambda_N}}^{\b,0}$.
From the monotonicity with respect to the volume and \cite{accn}, if $\beta<\beta_c$ 
there are constants $0<C(\beta)\leq C'(\beta)<\infty$
such that 
$
\left\langle
\s_i\s_j
\right\rangle_{\Lambda_N}^{\b,0}$ with free boundary condition satisfies
\[
\frac{C(\beta)}{|i-j|^{\alpha}}
\leq
\left\langle \s_i\s_j\right\rangle_{{\Lambda_N}}^{\b,0}
\leq
\left\langle \s_i\s_j\right\rangle_{{\mathbb{N}}}^{\b,0}
\leq 
\frac{C'(\beta)}{|i-j|^{\alpha}},
\]
where $C(\beta)\equiv (\beta \tanh\beta_c)/\beta_c$.
The uniformity of the lower bound is a simple application of FKG inequality. 
Using Lemma \ref{lem-esp} we have 
\[
N+\sum_{1\leq i<j\leq N}\frac{2C(\beta)}{|i-j|^{\alpha}}
\leq
\mathbb{E}_{\mathbb{P}_{N}^{2\b,0}}[\|S_N\|^2]
\leq
N+ 
\sum_{1\leq i<j\leq N}
\frac{2C'(\beta)}{|i-j|^{\alpha}}.
\]
Inequality (\ref{diffu}) follows by observing that $\sum_{1\leq i<j\leq N} \frac{1}{|i-j|^{\alpha}}={\cal O}(N)$.
\eop

\section{Central Limit Theorem}
\indent 

To derive a CLT for  \eqref{model} we make use of Theorem 1.2 and Theorem 3.1 from
\cite{WL}. It is required that a $C^2$-regularity condition to be satisfied. 
We say that a sequence of probability measures
 $\{\mu_N\}$ 
 satisfies the $C^2$-regularity condition if 
 for $Y_N$  
 with probability measure $\mu_N$  the following limit exists
\[
\Psi(t) =\ds\lim_{N\to\infty} \Psi_N(t)
          =\ds\lim_{N\to\infty} \frac{1}{N} \ln\mathbb{E}_{\mu_N}[\exp (tNY_N)].
\]
Moreover, for some neighborhood $[-\delta,\delta]$ of zero we have $\Psi (\cdot)<\infty$ and
 \[
 \Psi_N''(t) \underset{N\to\infty}{\longrightarrow} \Psi''(t) \  
 	\text{uniformly on  $[-\delta,\delta]$}. 
 \]
Under these hypotheses $Y_N$ is asymptotically Gaussian.  

For our polymer model take $v\in\mathbb{R}^2$ fixed and consider the empirical field projection
$$
L_N=\frac{1}{N} \langle S_N,v\rangle =\frac{1}{N}\sum_{j=1}^N \langle X_j,v\rangle. 
$$ 
Set $\mu_N = \mathbb{P}_N^{\beta,h}$ and define the pressure functional by
\[
 \Psi^{\b,h,v}(t) = \ds\lim_{N\to\infty} \Psi_N^{\b,h,v}(t) 
 						 =\ds\lim_{N\to\infty} \frac{1}{N} 
 							\ln \mathbb{E}_{\mathbb{P}_N^{\beta,h}}\left[\exp(\beta t N L_N)\right]
\]

\begin{theorem}\label{teo:clt}
Assume that the interactions are translation invariant and summable, 
that is, $V_{ij}=V(|i-j|)>0$ and $\sum_{i\in\mathbb{N}} V(i)<\infty$. For
$h\in\mathbb{R}^2$ with $h_1h_2\neq 0$ and any fixed  $v\in\mathbb{R}^2$ we have
$$
\frac{1}{\sqrt{N}}\left[ \b \langle S_N,v\rangle - N\mathbb{E}_{\mathbb{P}_N^{\beta,h}}[\b \langle S_N,v\rangle]\right]
\overset{\mathscr{D}}{\rightarrow} N\left(0,\dfrac{\partial^2}{\partial t^2}\Psi^{\b,h,v}(0)\right)
$$   
where  `` $\overset{\mathscr{D}}{\rightarrow}$ " stands for convergence in distribution.
\end{theorem}

\noindent {\bf Proof.} Under the hypotheses, the existence of the limit $\Psi^{\b,h,v}(\cdot)$ is proved in  \cite{R-TF}. To complete the $C^2$-regularity verification take complex number $z\in\mathbb{C}$ and express $\Psi_N^{\b,h,v}(z)$  
in terms of partition functions of one-dimensional Ising model.  
As in  Lemma \ref{lem-esp} write $\mathbb{Z}_N^{\b}(h) = Z_{\L_N}^{\b/2}(h_1)Z_{\L_N}^{\b/2}(h_2)$. 
Using the principal-value logarithm identities
$$
\ln(zw)= \ln z+\ln w+2\pi i\mathcal{K}(\ln z+\ln w)
$$
$$
\ln(z/w)= \ln z-\ln w+2\pi i\mathcal{K}(\ln z-\ln w)
$$
where $\mathcal{K}(x+iy)=-\sum_{n \geq -1} n I ( (2n-1)\pi < y \leq (2n+1)\pi)$ with $I(\cdot)$ being the indicator function, we have for  $v_1=\langle v, e_1-e_2\rangle$ and $v_2=\langle v, e_1+ e_2\rangle$
\begin{align*}
\Psi_N^{\b,h,v}(z)
=&\frac{1}{N}\left[\ln \mathbb{Z}_N^{\b}(h+zv)-\ln \mathbb{Z}_N^{\b}(h)
+2\pi i\mathcal{K}\left(\ln \mathbb{Z}_N^{\b}(h+zv)- \ln\mathbb{Z}_N^{\b}(h) \right)\right]
\\
=& \frac{1}{N}\ln{Z_{\L_N}^{\b/2}(h_1+zv_1)} + \frac{1}{N}\ln{Z_{\L_N}^{\b/2}(h_2+zv_2)} \nonumber 
\\
&+\frac{2\pi i}{N}\mathcal{K}\left( \ln Z_{\L_N}^{\b/2}(h_1+zv_1)+ \ln Z_{\L_N}^{\b/2}(h_2+zv_2)  \right) \nonumber
\\
& - \frac{1}{N}\ln{Z_{\L_N}^{\b/2}(h_1)} - \frac{1}{N}\ln{Z_{\L_N}^{\b/2}(h_2)} 
+\frac{2\pi i}{N}\mathcal{K}\left(\ln \mathbb{Z}_N^{\b}(h+zv)\right). 
\end{align*}
By assuming that
$ \text{Re}(h_i + zv_i) \neq 0$ and $h_1h_2 \neq 0$
 it follows from Lee-Yang Theorem's and standard  arguments from  \cite[p. 111]{R-SM} that 
$$
\Psi_N^{\b,h,v}(z) 
\to 
\Psi^{\b,h,v}(z),
\ \text{locally uniformly in}\ z.
$$  
Also, it follows that the derivatives of $\Psi_N^{\b,h,v}(z)$ 
converge uniformly on the compact subsets of $\mathbb{C}$. Hence the $C^2$-regularity condition is satisfied.
Since  
$$
\left.\frac{\partial^ 2}{\partial t^2} \Psi^{\beta,h,v}_N(t)\right|_{t=0}
=
\frac{1}{N}\  
			\mathbb{E}_{\mathbb{P}_N^{\beta,h}}\!\!
				\left[ \b\langle S_N,v\rangle - N\mathbb{E}_{\mathbb{P}_N^{\beta,h}}[\b\langle S_N,v\rangle] \right]^2 \to 
\frac{\partial^ 2}{\partial t^2} \Psi^{\beta,h,v}(0)
$$
we conclude the proof  using Theorem 3.1 from  \cite{WL}.
\eop

\begin{remark}
We emphasize that
in the above theorem we proved more than $C^2$-regularity condition.
In fact, we proved that the pressure is analytic. 
Another way to obtain the $C^2$-regularity condition for our polymer model  
is to apply both FKG and GHS inequalities, see \cite[p. 426]{WL}.   
\end{remark}

\section{Concluding Remarks}
\indent 

The random polymer model considered here 
interpolates between the SRW (infinite temperature) 
and a deterministic straight line (zero temperature). 
At very high temperatures 
this random polymer should be 
recurrent and transience would occur at very low temperatures, so we expect a recurrence-transience phase transition.
It would be interesting to prove the existence of 
such phase transition and also to determine the 
critical temperature that separates these two regimes.
\vspace*{2mm}

\noindent{\bf Acknowledgments.} We are grateful to L.R. Fontes for his many valuable comments and
careful reading of this manuscript.

\end{document}